%Paper: alg-geom/9506013
%From: Gregory Sankaran <G.K.Sankaran@pmms.cam.ac.uk>
%Date: Tue, 20 Jun 95 15:28:33 BST

\input mssymb
\def\CC {{\Bbb C}}
\def\FF {{\Bbb F}}
\def\HH {{\Bbb H}}

\def\PP {{\Bbb P}}
\def\QQ {{\Bbb Q}}
\def\RR {{\Bbb R}}
\def\ZZ {{\Bbb Z}}
\def\contin {\subseteq}

\def\Bar{\overline}

\def\minus {\setminus}
\def\cross {\times}
\def\part#1#2 {{\partial {#1}/\partial {#2}}}

\def\Sing {\mathop{\rm Sing}\nolimits}
\def\Aut {\mathop{\rm Aut}\nolimits}

\def\Hom {\mathop{\rm Hom}\nolimits}

\def\Sp {\mathop{\rm Sp}\nolimits}
\def\SL {\mathop{\rm SL}\nolimits}
\def\PSL {\mathop{\rm PSL}\nolimits}
\def\GL {\mathop{\rm GL}\nolimits}
\def\Ker {\mathop{\rm Ker}\nolimits}

\def\diag {\mathop{\rm diag}\nolimits}
\def\codim {\mathop{\rm codim}\nolimits}

\def\im {\mathop{\rm Im}\nolimits}

\def\cH {{\cal H}}

\def\bdn {{\bf n}}

\def\bdx {{\bf x}}
\def\bdy {{\bf y}}

\def\into {\hookrightarrow}
\def\To{\longrightarrow}

\def\Prod{\prod\limits}

\def\tens{\otimes}

\def\remark{\noindent{\sl Remark. }}
\def\pf{\noindent{\sl Proof: }}

\font\gothic=eufm10
\def \gothg {{\hbox{\gothic g}}}

\def \hf {{{1}\over{2}}}
\def\rationalmap{\mathrel{{\hbox{\kern2pt\vrule height2.45pt depth-2.15pt
 width2pt}\kern1pt {\vrule height2.45pt depth-2.15pt width2pt}
  \kern1pt{\vrule height2.45pt depth-2.15pt width1.7pt\kern-1.7pt}
   {\raise1.4pt\hbox{$\scriptscriptstyle\succ$}}\kern1pt}}}
\def\qed{\vrule width5pt height5pt depth0pt\par\smallskip}
\def\surj{\to\kern-8pt\to}
\outer\def\startsection#1\par{\vskip0pt
 plus.3\vsize\penalty-100\vskip0pt
  plus-.3\vsize\bigskip\vskip\parskip\message{#1}
   \leftline{\bf#1}\nobreak\smallskip\noindent}
\outer\def\camstartsection#1\par{\vskip0pt
 plus.3\vsize\penalty-150\vskip0pt
  plus-.3\vsize\bigskip\vskip\parskip
   \centerline{\it#1}\nobreak\smallskip\noindent}
\def\chain{\dot{\hbox{\kern0.3em}}}
\def\cochain{\d{\hbox{\kern0.3em}}}

\def\imic{\cong}
\outer\def\thm #1 #2\par{\medbreak
  \noindent{\bf Theorem~#1.\enspace}{\sl#2}\par
   \ifdim\lastskip<\medskipamount \removelastskip\penalty55\medskip\fi}
\outer\def\prop #1 #2\par{\medbreak
  \noindent{\bf Proposition~#1.\enspace}{\sl#2}\par
   \ifdim\lastskip<\medskipamount \removelastskip\penalty55\medskip\fi}
\outer\def\lemma #1 #2\par{\medbreak
  \noindent{\bf Lemma~#1.\enspace}{\sl#2}\par
   \ifdim\lastskip<\medskipamount \removelastskip\penalty55\medskip\fi}
\outer\def\coro #1 #2\par{\medbreak
  \noindent{\bf Corollary~#1.\enspace}{\sl#2}\par
   \ifdim\lastskip<\medskipamount \removelastskip\penalty55\medskip\fi}
\def\deep #1 {_{\lower5pt\hbox{$#1$}}}

\baselineskip=18pt

\raggedbottom
\magnification=1200

\def \Ga {{\Gamma}}
\def \onto {{\surj}}
\def \UZ {{U(F)_{\ZZ}}}
\def \DB {{\Bar{D/\UZ}}}
\def \GM {{\tilde\Ga^0_{1,3}(2)}}

\centerline{\bf Fundamental Group of Locally Symmetric Varieties}
\medskip
\centerline{\it G.K. Sankaran}
\bigskip
\noindent The geometry of moduli spaces of complex abelian varieties
and of compactifications of those moduli spaces has been the object of
much study in the last few years. Given such a compactified moduli
space it is natural to ask about the fundamental group of a resolution
of singularities. This problem has been studied in some special cases,
for instance in [K], [HK] and [HS]. It follows from the results of
[K] and the well-known fact that every arithmetic subgroup of
$\Sp(2g,\QQ)$ is a congruence subgroup of some level that the
fundamental group must be finite except when~$g=1$, i.e., except in
the case of modular curves. The method in all these cases, and here,
is to use the toroidal compactification of [SC], considering the
moduli spaces as quotients of the Siegel upper half-space by
arithmetic subgroups of the symplectic group.

In this paper we treat the subject in greater generality. In many
cases we are able to identify the fundamental group explicitly as a
quotient of the arithmetic group in question. For this purpose we do
not need to restrict ourselves to the symplectic group but instead may
consider any locally symmetric variety. Later we return to the case of
moduli of abelian varieties (specifically, to Siegel modular
varieties) and calculate the fundamental group in some interesting
special cases. In many cases, including those studied in the papers
mentioned above, the fundamental group is trivial, but we give
examples to show that this need not be true in general.

The example in Proposition~3.1 is a modified version of one suggested
to me by Professor M.S.~Raghunathan. I am grateful to him for pointing
it out to me. I am also grateful to Professors K.~Hulek and W.~Ebeling
for useful remarks, and to T\^ohoku University and especially to
Professor Tadao Oda for their hospitality during a visit to Japan when
I began this work.

\startsection 1. Unipotent elements in parabolic subgroups

We shall study locally symmetric varieties and their compactifications
as described in~[SC]. We therefore let $D$ be a bounded symmetric
domain with $\Aut (D)^0=G$, a simple real Lie group defined
over~$\QQ$. We fix an arithmetic subgroup $\Ga$ of~$G$. We put
$X=D/\Ga$ and use the methods of [SC] to construct a
compactification $\bar X$ of~$X$. Since we have not required $\Ga$
to be a neat subgroup of~$G$ we cannot guarantee that we can
choose $\bar X$ to be smooth, but we can take a resolution of
singularities $\tilde X\to\bar X$ if we wish. The object of interest to
us is the topological fundamental group~$\pi_1(\tilde X)$. Some
arbitrary choices have to be made in constructing $\bar X$ and $\tilde
X$ but $\pi_1(\tilde X)$ does not depend on the choices.

We begin as in [HS] with some topological lemmas, repeated
here for ease of reference.

\lemma 1.1 Let $M$ be a connected, simply-connected real manifold and
$G$ a group acting discontinuously on~$M$. Take a base point $x\in M$.
Then the quotient map $\phi:M\to M/G$ induces a surjective
homomorphism $\phi_*:G\to\pi_1\big(M/G,\phi(x)\big)$.

\pf [G], Satz~5. The quotient map~$\phi_*$ is easy to describe
and this is done in~[HS].~\qed

\lemma 1.2 Let $M$ be a connected complex manifold and $M_0$ an
analytic subvariety. Let $x\in M\minus M_0$. Then the inclusion
$M\minus M_0\into M$ induces a map $\pi_1(M\minus
M_0,x)\to\pi_1(M,x)$, which is surjective if $\codim_\CC M_0\ge 1$ and
an isomorphism if $\codim_\CC M_0\ge 2$.

\pf See [HK].~\qed

\lemma 1.3 If $X_1$ and $X_2$ are bimeromorphically equivalent
connected complex manifolds, then $\pi_1(X_1)\imic\pi_1(X_2)$.

\pf Also in [HK].~\qed

Now let $\phi:D\to D/\Ga=X$ be the quotient map and let $X_{\rm reg}$
be the nonsingular locus of~$X$. Put $D'=\phi^{-1}(X_{\rm reg})$. The
complement $D\minus D'$ consists of countably many analytic
subspaces of~$D$ of complex codimension at least~$2$; recall that $D$
has naturally the structure of a complex manifold. So $D'$ is
simply-connected and we get a surjection $\Ga\onto\pi_1(X_{\rm reg})$.
Composing this with the surjection $\pi_1(X_{\rm reg})\onto\pi_1(\tilde
X)$ induced by the inclusion $X_{\rm reg}\into\tilde X$ we get a
surjection $\psi:\Ga\onto\pi_1(\tilde X)$.

\thm 1.4 If $P$ is a parabolic subgroup of $G$ and $U_P$ is the centre
of the unipotent radical of~$P$ then $U_P\cap\Ga\subseteq\Ker\psi$.

\pf We can assume that $P\cap\Ga\ne 1$. By [SC], III.3.2,
Proposition~2, $P$~corresponds to a boundary component~$F_P=F$.
Since $P\cap\Ga\ne 1$, $F$~is a rational boundary component and we may
adopt the notation of~[SC] and write $U_P\cap\Ga=\UZ$.

Recall the procedure of partial compactification in the direction
of~$F$ as it is described in~[SC]: we can write
$$
D={\rm Im}^{-1}C(F)\cross\CC^k\cross F\contin
D(F)=U(F)_\CC\cross\CC^k\cross F,
$$
where $U(F)=U_P$, $U(F)_\CC$ is the complexification of $U(F)$,
$C(F)$ is an open convex cone in $U(F)$ and ${\rm Im}:U(F)_\CC\to U(F)$
is the imaginary part. $\UZ$ acts by translation in the real direction
in $U(F)_\CC$ and hence preserves~$C(F)$: it acts trivially on
$\CC^k\cross F$ and there is an embedding
$$
D/\UZ\into T(F)\cross\CC^k\cross F
$$
where $T(F)=U(F)_\CC/\UZ$ is an algebraic torus over~$\CC$. By
choosing an appropriate fan $\{\sigma_\alpha\}$ subdividing~$C(F)$ we
construct a partial compactification $T(F)_{\{\sigma_\alpha\}}$
of~$T(F)$ and hence
$$
\DB=\Big(D/\UZ\Big)_{\{\sigma_\alpha\}} \contin
T(F)_{\{\sigma_\alpha\}}\cross\CC^k\cross F.
$$
We can even do this in such a way as to make $\DB$ smooth and keep
everything $\Ga$-equivariant.

$T(F)_{\{\sigma_\alpha\}}$ is simply-connected: see for instance [F].
If $T(F)=\Hom (M,\CC^*)$, where $M$ is a lattice, and if $N$ is the
dual lattice, then $\pi_1\big(T(F)\big)$ consists of the classes of
the loops $s\mapsto\exp\{2\pi is<\ ,\bdn >\}$ for $\bdn\in N$. If
$\bdn\in N\cap\sigma_\alpha$ then this loop is killed in
$\pi_1\big(T(F)_{\{\sigma_\alpha\}}\big)$ by a retraction given by
$$
R_{\bdn}(s,t)=t\exp\{2\pi is<\ ,\bdn>\}\qquad
\hbox{for }(s,t)\in[0,1]^2,t\not=0,
$$
extended by putting $R_\bdn(s,0)=\lim\limits_{t\to 0}R_\bdn(s,t)$, which
exists in $T(F)_{\{\sigma_\alpha\}}$ (and is independent of~$s$).
Since $N\cap\sigma_\alpha$ generates $N$ if $\sigma_\alpha$ is of
maximal dimension, $T(F)_{\{\sigma_\alpha\}}$ is simply-connected.

We can identify $T(F)$ with $N\tens_\ZZ\CC^*$, so that
$R_\bdn(s,t)=\bdn{\scriptstyle\tens}te^{2\pi is}$. We can consistently
choose a logarithm for $t\in(0,1]$, so that $R_\bdn$ comes from a map
$\hat R_\bdn:[0,1]\cross(0,1]\to N\tens_\ZZ\CC$ given by $\hat
R_\bdn(s,t)=\bdn{\scriptstyle\tens}(\log t +2\pi is)$. Note also that
$\im\hat R_\bdn(s,t)\in C(F)$ if $\bdn\in C(F)$. From this it follows
that $\DB$ is also simply-connected, since $F$ and $\CC^k$ are both
simply-connected themselves. As $\DB$ maps under the action of~$\Ga$
onto a Zariski-open set in $\bar X$, and as $\UZ$ acts trivially
on~$\DB$, we see that $\UZ$ is in the kernel of $\Ga\to\pi_1(\bar X)$.

In order to show that if $\eta\in\UZ$ then
$\psi(\eta)=1\in\pi_1(\tilde X)$ we have to do a little more: we must
avoid the singularities of~$\bar X$. Suppose $H:[0,1]^2\to
\DB$ is a null homotopy for a loop in~$\DB$ coming from
$\eta\in\UZ$, and let $\bar H$ be the corresponding null homotopy
in~$\bar X$. We may assume that~$H$ and~$\bar H$ are smooth maps. The
singularities of~$\bar X$ are quotient singularities arising from the
action of $\Ga$ on $\DB$, since we have chosen $\{\sigma_\alpha\}$ so
as not to lead to singularities in~$\DB$. We may assume, by the
description of $R$ above, that there is a unique point $x_0$ (which we
can take to be $(\hf,\hf)$) in $[0,1]^2$ such that
$H(x_0)\not\in D/\UZ$, and that away from this point $H$~lifts from a map
into $D/\UZ$ to a map
$$
\hat H:[0,1]^2\minus\{({\scriptstyle{\hf},{\hf}})\}\To D.
$$

Denote by $Z$ the preimage in $\DB$ of $\Sing \bar X$. Then~$Z$ is the
union of countably many analytic submanifolds of codimension~$\ge 2$.
The same is true of~$\hat Z$, the preimage of $Z\cap\big(D/\UZ\big)$
in~$D$. By choosing a path from~$x$ to~$\eta(x)$ in $D\minus\hat Z$
we may assume that $\hat H(0,t)$ misses~$\hat Z$ and that the
loop~$H(0,t)$ misses~$Z$. The idea now is to move~$H$ so that its
image misses~$Z$, without changing~$H$ on $\{0,1\}\cross\{0,1\}$: this
will produce a null homotopy of the path $\bar H(0,t)$ corresponding
to~$\eta$ that takes place entirely in $\bar X\minus\Sing\bar X$
and therefore lifts to~$\tilde X$.

To move the part of~$H$ that lies inside $D/\UZ$, we choose a
neighbourhood $N_1$ of $1\in G$ and a diffeomorphism
$$
\theta_1:N_1\To\Delta_1=\{\bdx\mid |\bdx|<1\}\contin {\gothg}={\rm Lie}\, G.
$$
Then we define $\hat\cH_1:[0,1]^2\cross\Delta_1\to D$ by
$$
\hat\cH_1\big((a,b),\bdx\big)=\theta_1^{-1}\big(\lambda_1(a,b)\big)\hat
H(a,b),
$$
where $\lambda_1:\RR^2\to\RR$ is a $C^\infty$~function satisfying
$0\le\lambda(a,b)<1$ and $\lambda_1(a,b)=0$ if and only if $a=b=\hf$
or $(a,b)\not\in[0,1]^2$. Let
$\hat\cH_1:[0,1]^2\cross\Delta_1\to D/\UZ$ be the composition
with the quotient map.

To move~$H$ near the boundary, choose a neighbourhood $N_\infty$ of
$H\big((\hf,\hf)\big)$ and a diffeomorphism
$$
\theta_\infty:N_\infty\To\Delta_\infty=\{\bdx\mid|\bdx|<1\}
\contin\RR^{2\dim_\CC D}.
$$
Let $\Delta$ be an open disc centred at $(\hf,\hf)$ such that
$H(\Delta)\contin N_\infty$ and $\theta_\infty
H(\Delta)\contin\hf\Delta_\infty$. We translate~$H$ near~$(\hf,\hf)$,
defining $\cH_\infty:[0,1]^2\cross\Delta_\infty\to\DB$ by
$$
\cH_\infty\big((a,b),\bdx\big)=\theta_\infty^{-1} \big(\theta_\infty
H(a,b)+\lambda_\infty(a,b)\bdx\big),
$$
where $\lambda_\infty:\RR^2\to\RR$ is a $C^\infty$~function satisfying
$0\le\lambda_\infty(a,b)<\hf$ and $\lambda_\infty(a,b)=0$ if and only
if $(a,b)\not\in\Delta$.

Now we compose $\cH_1$ and $\cH_\infty$ so as to get a map
$\cH:[0,1]^2\cross\Delta_1\cross\Delta_\infty\to\DB$ given by
$$
\cH\big((a,b),\bdx,\bdy\big)=\theta_\infty^{-1} \Big\{\theta_\infty
\cH_1\big((a,b),\bdx\big)+\lambda_\infty(a,b)\bdy\Big\}.
$$
Since the real tangent space to~$D/\UZ$ at any point is a quotient
of~$\gothg$ as a real vector space, $\cH$~is a submersion for
$(a,b)\in(0,1)^2\minus\{(\hf,\hf)\}$. By the choice of $\cH_\infty$
it is a submersion at $(\hf,\hf)$ as well. So by the Thom
transversality theorem (see [BG]) there exist $\bdx_0\in\Delta_1$ and
$\bdy_0\in\Delta_\infty$ such that $\cH\big((s,t),\bdx_0,\bdy_0\big)$
is a null homotopy tranversal to each component of~$Z$. As~$Z$ has real
codimension~$4$ this means the homotopy misses~$Z$, and we are
done.~\qed

\remark The technical difficulties above arose almost entirely from
the possibility that $\bar X$ might have singularities at or near the
boundary. If~$\Ga$ is neat, so that $\bar X$ is smooth, or if, as in
[HS], we need only consider boundary components far from $\Sing\bar
X$, the situation is very simple. Even if we cannot easily avoid
$\Sing\bar X$, it is often the case that the resolution $\tilde
X\to\bar X$ can be chosen to have simply-connected fibres (for
instance if $\bar X$ has only cyclic quotient singularities), and then
$\pi_1(\tilde X)\imic\pi_1(\bar X)$ anyway. Perhaps this can always be
done, even when $\bar X$ has unknown finite quotient singularities.

For practical purposes, the following easy corollary, or something
like it, is often useful.

\coro 1.5 The fundamental group $\pi_1(\tilde X)$ is a quotient of
$\Pi(\Gamma)=\Ga/\Upsilon$, where $\Upsilon$ is the normal subgroup of
$\Ga$ generated by all elements in the centre of the unipotent radical
of some parabolic subgroup of~$G$.

\pf $\Pi(\Ga)$ is the largest quotient of $\Ga$ satisfying the
conditions of Theorem~1.4.~\qed

In fact if we let $\tilde\Ga$ be the normaliser of $\Ga$ in $G(\QQ)$
then $\tilde\Ga/\Ga$ acts on~$X$ and preserves the set of rational
boundary components; the action therefore extends to $\bar X$ if we
choose the compactification correctly. The resolution $\tilde X\to\bar
X$ may also be chosen to be $\tilde\Ga/\Ga$-equivariant, and therefore
$\tilde\Ga/\Ga$ acts on $\pi_1(\tilde X)$. So the kernel of
$\psi:\Ga\onto\pi_1(\tilde X)$ must be a normal subgroup
of~$\tilde\Ga$. But in fact $\Pi(\Ga)$ satisfies this condition as
well. For if $g\in\tilde\Ga$ and $P$~is a parabolic subgroup of~$G$
defined over~$\QQ$, then so is~$P^g$ and the unipotent radical $R_u$
satisfies $R_u(P^g)=R_u(P)^g$. So $U_{P^g}=U_P^g$. This shows also
that $\Upsilon$ is actually the subgroup of~$\Ga$ generated by all
the~$U_P$s, which is automatically normal. But for purposes of
calculation the statement in Corollary~1.5 is preferable, because
rather than work with all parabolic subgroups it is easier to take a
representative from each conjugacy class. Geometrically, this means
that one works with the boundary components of $\bar X$ (i.e., with
the Tits building) rather than with boundary components of~$D$. It is
not even necessary to look at all boundary components: in view of the
fact ([SC], III.4.4, Theorem~3) that $U(F)\contin U(F')$ if $F'$ is a
boundary component of~$F$ it is enough to look at the minimal ones.
But it is usually easiest to write down $U(F)$ when $F$ is a maximal
proper boundary component, corresponding to a divisor in the toroidal
compactification. This is often sufficient because it happens that
such $U(F)$ already generate $U(F')$ for all~$F'$.

Finally, the condition that $G$ is simple is also stronger than we
need. If we allow $G$ to be semisimple then nothing changes except
that we have (cf.~[SC], p.208) $D=\Prod D_i$ (for instance $D=\HH^2$
and $X$ is a Hilbert modular surface) and
$G=\Prod\Aut(D_i)^0=\Prod G_i$ (as a set), and we must allow $P=\Prod
P_i$, where each $P_i$ is either parabolic or else equal to~$G_i$.

We summarise our results as follows.

\coro 1.6 Let $D$ be a bounded symmetric domain with $G=\Aut(D)^0$ a
semisimple real Lie group defined over~$\QQ$. Let $\Gamma$ be an
arithmetic subgroup of~$G$. Then $\pi_1(\tilde X)$ is a quotient of
$\Pi(\Ga)=\Ga/\Upsilon$, where $\Upsilon$ is the subgroup generated by
$U(F)\cap\Ga$ for all rational boundary components~$F$ of~$D$.
Equivalently, $\Upsilon$ is the normal subgroup of~$\Ga$ generated by
$U(F)\cap\Ga$ as~$F$ runs through a set of representatives for all
boundary components of~$\bar X$ (or all minimal boundary components
of~$\bar X$).~\qed

\startsection 2. Neat arithmetic groups.

Throughout this section we assume that~$\Ga$ is neat. In particular
this implies that $\Ga$ is torsion-free. In this special (but not very
special) case we can identify $\pi_1(\tilde X)$ precisely. Note that
now~$\tilde X=\bar X$.

\thm 2.1 Let $D$ be a bounded symmetric domain and $G=\Aut(D)^0$ a
semisimple Lie group defined over~$\QQ$. Let $\Ga$ be a neat
arithmetic subgroup of~$G$. Then $\pi_1(\tilde X)\imic\Pi(\Ga)$.

\pf Let $\Upsilon$ be as in Corollary~1.6 above. We should like to
proceed by putting $X_u=D/\Upsilon$ and constructing a toroidal
compactification~$\bar X_u$ (which could be assumed to be smooth
because $\Upsilon$, being a subgroup of~$\Ga$, is neat). Then we should
expect $\Ga/\Upsilon$ to act freely on $\bar X_u$ with quotient $\bar
X$. Unfortunately $\Upsilon$ need not have finite index in~$\Ga$ so
$\Upsilon$ will not, in general, be an arithmetic group. But this is
not a serious difficulty: the only consequence is that $\bar X_u$ need
not be compact, which does not matter to us.

We construct $\bar X_u$ as in [SC],~III.5. For each rational boundary
component~$F$ we have $\Ga\cap U(F)=\Upsilon\cap U(F)$, by the
definition of~$\Upsilon$, so we can call this group~$\UZ$ without
ambiguity. We take the same fans $\{\sigma_\alpha\}$ as before,
requiring them to be admissible for the action of~$\Ga$, not just
of~$\Upsilon$, and we define $\bar X_u$ as the quotient of
$\coprod\limits_F\DB$ by the closure of the equivalence relation defined by
the action of~$\Upsilon$. Thus (cf.~[SC], p.255) if
$x_i\in\Bar{D/U(F_i)_\ZZ}$, $i=1,$~$2$, then $x_1\sim x_2$ if and only
if there are a rational boundary component~$F$, an element
$\eta\in\Upsilon$ and a point $x\in\DB$ such that $F_1$ and $\eta F_2$
are boundary components of~$F$ and $x$ projects to $x_1$ and to $\eta
x_2$ under the projection maps $\DB\to\Bar{D/U(F_1)_\ZZ}$ an
d$\DB\to\Bar{D/U(\eta F_2)_\ZZ}$. This equivalence relation is closed,
so $\bar X_u$ is Hausdorff: indeed, $\pi_F^{(u)}:\DB\to\bar X_u$ is
biholomorphic onto its image. But the map $\pi_F:\DB\to\bar X$ factors
as $\pi_F=q\pi_F^{(u)}$, where $q:\bar X_u\to\bar X$ is the quotient
map under the action of~$\Ga/\Upsilon$, and $\pi_F$ is \'etale since
$\Ga$ is neat. Therefore $q$~is \'etale and so
$\Ga/\Upsilon\into\pi_1(\bar X)$.~\qed

\startsection 3. Examples.

If $X$ is a curve then $\pi_1(\tilde X)$ will not be finite unless
$\tilde X=\PP^1$. In this case (for instance when $G=\SL(2,\RR)$ and
$\Ga$ is the principal congruence subgroup of some level $l\ge 5$)
$\Upsilon$ is of infinite index in~$\Ga$ and in particular is not an
arithmetic group. In higher-dimensional cases (other than products
with a factor of this type) $\Ga/\Upsilon$ seems to be finite, but I
am not aware of any definite general result of that nature. In many
cases $\Ga/\Upsilon$ is trivial and thus $\tilde X$ is
simply-connected: this is shown for various Siegel modular varieties
in [K], [HK] and [HS]. The most frequently considered Hilbert modular
surfaces are also simply-connected but some others are not: see~[vdG].
Likewise, Siegel modular threefolds are not simply-connected in
general, and the fundamental group can even be quite big, as the
following examples show.

\thm 3.1 Let $l\ge 4$ be an integer and let $p$~be a prime not
dividing~$2l$. take $\Ga(l)$ to be the principal congruence subgroup
of level~$l$ and let $\Ga(l)_p\contin\Sp(4,\FF_p)$ be the image of
$\Ga(l)$ under the reduction mod~$p$ map ${\rm
red}_p:\Ga(l)\to\Sp(4,\FF_p)$. Let $\Ga^q(l)_p$ be a Sylow
$q$-subgroup of $\Ga(l)_p$ for some prime $q\ne p$, and let $\Ga={\rm
red}_p^{-1}\big(\Ga(l)_p\big)$. Then if $X=D/\Ga$ (in this case $D$ is
the Siegel upper half-plane), $\pi_1(\tilde X)$ has a quotient
isomorphic to $\Ga^q(l)_p$.

\pf We must check first that the groups mentioned exist, by showing
that $\Ga(l)_p$ is not a $p$-group. The map
$$
{\rm red}_p:\Sp(4,\ZZ)\To\Sp(4,\FF_p)
$$
is surjective, so given $\alpha\in\Sp(4,\FF_p)$ and $\beta\in\Ga(l)_p$
we can choose $\tilde\alpha\in\Sp(4,\ZZ)$ and $\tilde\beta\in\Ga(l)$
such that ${\rm red}_p(\tilde\alpha)=\alpha$ and ${\rm
red}_p(\tilde\beta)=\beta$. By definition $\Ga(l)$ is a normal
subgroup of $\Sp(4,\ZZ)$, so
$\tilde\alpha^{-1}\tilde\beta\tilde\alpha\in\Ga(l)$ and therefore
$\alpha^{-1}\beta\alpha\in\Ga(l)_p$. So $\Ga(l)_p$ is a normal
subgroup of $\Sp(4,\FF_p)$. But the only nontrivial normal subgroup of
$\Sp(4,\FF_p)$ is its centre, which has order~$2$, and $\Ga(l)_p$
contains the element $\pmatrix{I&{\rm red}_p(l)I\cr 0&I\cr}$, which
has order~$p\ne 2$. So in fact $\Ga(l)_p$ is the whole of $\Sp(4,\FF)$
and the order of this group is not a power of~$p$.

$\Ga(l)$ is neat and so, therefore, is~$\Ga$. Moreover, $\Ga$ is
evidently of finite index in~$\Ga(l)$ and is thus an arithmetic group.
So by Theorem~2.1 we can construct a smooth compactification $\bar X$
of~$X$ whose fundamental group is $\Ga/\Upsilon$, where $\Upsilon$ is
a group generated by unipotent elements. Suppose $\tilde\eta$ is a
unipotent element of $\Ga$. Then ${\rm red}_p(\tilde\eta)=\eta$ is a
unipotent element of $\Sp(4,\FF_p)$ and is therefore conjugate over
$\GL(4,\FF_p)$ to an upper-triangular element $\eta'\in\GL(4,\FF_p)$.
The order of~$\eta'$ is obviously a power of~$p$, so the order
of~$\eta$ is a power of~$p$, but $\eta\in\Ga^q(l)_p$ which is a
$q$-group, so $\eta$ is the identity. So $\Upsilon\contin\Ker{\rm
red}_p$ and ${\rm red}_p$ gives a surjective morphism
$\Ga/\Upsilon\onto\Ga^q(l)_p$.~\qed

\noindent{\sl Remarks} i) We can even take $p=2$ if we like, since
$\pmatrix{I&{\rm red}_2(l)I\cr 0&I\cr}$ is not central.

\noindent\phantom{\it Remarks} ii) The varieties exhibited in
Theorem~3.1 are all of general type because of the result of Yamazaki
([Y]) that Siegel modular threefolds of level~$l$ are of general type
for~$l\ge 4$.

\coro 3.2 The fundamental group of a Siegel modular threefold need not
be abelian.

\pf $\Ga^q(l)_p$ can be any $q$-subgroup, not necessarily a Sylow
$q$-subgroup. If we take $l=5$, $p=7$ and $q=2$ then we
can take $\Ga^q(l)_p$ to be non-abelian, since $\Sp(4,\FF_7)$ contains
a subgroup isomorphic to $\PSL(2,\FF_7)\imic\PSL(4,\FF_3)$, which
obviously contains a dihedral group of order~$8$.~\qed

In spite of the possibility that $\pi_1(\tilde X)$ may be a large
finite group, it is frequently the case that Siegel modular threefolds
arising in geometry turn out to be simply-connected. This is shown for
the principal congruence subgroups (in terms of abelian surfaces, the
moduli spaces of principally polarised abelian surfaces with level~$l$
structure) in~[K] and for the moduli of $(1,p)$-polarised abelian
surfaces with level structure ($p$~an odd prime) in [HS]. Here we will
give two more cases where $\tilde X$ is simply-connected and one
where it is not.

\thm 3.3 Let $X$ be the moduli space of abelian surfaces with a
polarisation of type $(1,p)$, $p$~an odd prime. Then $\tilde X$ is
simply-connected.

\pf This is a fairly straightforward consequence of the calculations
in [HS] (or can be proved directly in the same way, using
Theorem~1.4). All we need to do is observe that, adopting the notations
of [HKW] and [HS], the kernel of $\psi:\Ga^0_{1,p}\to\pi_1(\tilde X)$ contains
not
only $M_0$ and $\Ga(p^2)$ but also the extra element
$$
M_0'=\pmatrix{ 1&0&0&0\cr
               0&1&0&p\cr
               0&0&1&0\cr
               0&0&0&1\cr}
    =j_2\pmatrix{1&1\cr 0&1\cr}.
$$
In [HS] it is shown that a normal subgroup of $\Ga^0_{1,p}$ containing
$M_0$ and $\Ga(p^2)$ must contain $\Ga_{1,p}$. The normal subgroup
$\Ker\psi$ is larger than that since $M_0'\not\in\Ga_{1,p}$, and must
therefore be equal to $\Ga^0_{1,p}$ since $\Ga^0_{1,p}/\Ga_{1,p}$ is
isomorphic to $\SL(2,\FF_p)$ and hence simple modulo~$\pm I$, which acts
trivially.~\qed

Two Calabi-Yau threefolds occur in the paper [BN] of Barth and Nieto.
One, called $N$ there and also described by Naruki in [Na], is a
certain singular quintic hypersurface in $\PP^4$; the other, called
$\tilde N$, is a double cover of~$N$ and is birationally equivalent to
the moduli space of abelian surfaces with a polarisation of type
$(1,3)$ (or $(2,6)$) and level-$2$ structure.

\thm 3.4 If $Z$ is a desingularisation of $\tilde N$ then
$\pi_1(Z)\imic\ZZ/2\cross\ZZ/2$.

\pf Put $E=\pmatrix{1&0\cr 0&3\cr}$ and $\Lambda=\pmatrix{0&E\cr
-E&0}$. The moduli space of $(1,3)$-polarised abelian surfaces with
level-$2$ structure is $\HH_2/\GM$, where
$\GM$ is the kernel of reduction mod~$2$ in
$\tilde\Ga^0_{1,3}=\Sp(\Lambda,\ZZ)$ and $\GM$ acts on~$\HH_2$ by
$$
\pmatrix{A&B\cr C&D\cr}:Z\To (AZ+B)(CZ+D)^{-1}E.
$$
One can check, by the same method and with the same notation as in
[HS], that $\GM$ is generated by $\tilde\jmath_1\big(\Ga_1(2)\big)$,
$\tilde\jmath_2\big(\Ga_1(2)\big)$, $\tilde M_1^2$, $\tilde M_2^2$,
$\tilde M_3^2$ and $\tilde M_4^2$. We must describe the normal
subgroup generated by the centres of unipotent radicals of parabolic
groups in terms of these. $\tilde M_0^2=\tilde\jmath_1\pmatrix{1&2\cr
0&1}$ is in the centre of the unipotent radical of the parabolic
subgroup corresponding to the $\Lambda$-isotropic line $\QQ(1,0,0,0)$,
and the normal subgroup of $\Ga_1(2)$ generated by $\pmatrix{1&2\cr
0&1}$ is the whole of $\Ga_1(2)$ (this reflects the fact that the
modular curve of level~$2$ is rational). So we get all of
$\tilde\jmath_1\big(\Ga_1(2)\big)$, and similarly we get
$\tilde\jmath_2\big(\Ga_1(2)\big)$ from ${\tilde M'_0}{}^2$ corresponding
to $(0,1,0,0)$. The elements $\tilde M_1^2$ and $\tilde M_2^2$ are
also unipotent and they lie in the centres of the unipotent radicals
of the parabolic groups corresponding to the isotropic planes spanned
by $(1,0,0,0)$ and $(0,1,0,0)$ and by $(0,0,1,0)$ and $(0,0,0,1)$
respectively. All other unipotent elements are conjugate to products
of these. We can also generate
$$
\tilde M_3^4=\tilde\jmath_2\pmatrix{\phantom{-}1&0\cr -2p&1}
\tilde M_0^2\tilde M_1^2\tilde M_0^{-2}\tilde M_1^{-2},
$$
and $\tilde M_4^4$ similarly, but not $\tilde M_3^2$ or $\tilde
M_4^2$. This is because
$$
\GM=\left\{\gamma\in\Sp(\Lambda,\ZZ)\mid\gamma-I\in\pmatrix{
          2\ZZ&2\ZZ&2\ZZ&2\ZZ\cr
          6\ZZ&2\ZZ&6\ZZ&2\ZZ\cr
          2\ZZ&2\ZZ&2\ZZ&2\ZZ\cr
          6\ZZ&2\ZZ&6\ZZ&2\ZZ\cr}
\right\}
$$
contains the normal subgroup
$$
\tilde\Upsilon^0_{1,3}(2)=
\left\{\gamma\in\Sp(\Lambda,\ZZ)\mid\gamma-I\in\pmatrix{
          \phantom{1}2\ZZ&4\ZZ&\phantom{1}2\ZZ&2\ZZ\cr
                    12\ZZ&2\ZZ&\phantom{1}6\ZZ&2\ZZ\cr
          \phantom{1}2\ZZ&2\ZZ&\phantom{1}2\ZZ&4\ZZ\cr
          \phantom{1}6\ZZ&2\ZZ&          12\ZZ&2\ZZ\cr}
\right\}.
$$
So in this case $\Pi\big(\GM\big)\imic\GM/\tilde\Upsilon^0_{1,3}(2)
\imic\ZZ/2\cross\ZZ/2$.

We cannot at once conclude that $\pi_1(Z)\imic\ZZ/2\cross\ZZ/2$
because $\GM$ is not neat. The elements that cause it to fail to
be neat are the conjugates of $-I$ (which acts trivially on $\HH_2$)
and of $\diag(1,-1,1,-1)=I_0$. But $I_0$ fixes a surface and acts,
even at the cusps, as a reflection at its fixed points, so that $\bar
X'$, the compactification of $\HH_2/\tilde\Upsilon^0_{1,3}(2)$, is
smooth even though $\tilde\Upsilon^0_{1,3}(2)$ is not neat.
Furthermore, $\Pi\big(\GM\big)$ acts freely on $\bar X'$ because if
the image in $\Pi\big(\GM\big)$ of $\tilde M_3^2$ (or $\tilde M_4^2$)
preserves some boundary component~$F$ and fixes $x\in\DB$ then
$\tilde M_3^2$ is in the group generated by $\UZ$, $I_0$ and $-I$, but
these are all in $\tilde\Upsilon^0_{1,3}(2)$. But $\bar X'$ is
simply-connected by construction, so $\pi(Z)\imic\ZZ/2\cross\ZZ/2$.~\qed

To calculate the fundamental group of a resolution of~$N$ we need a
different method, using the projective description of the variety.

\prop 3.4 Any desingularisation~$Y$ of~$N$ is simply-connected.

\pf $N$ is given by the equations
$$
u_0+u_1+u_2+u_3+u_4+u_5=\sum_{i=0}^5\prod_{j\ne i}u_j=0
$$
in~$\PP^5$. Taking the hyperplanes $H_0=(u_0=0)$ and $H_1=(u_1=0)$
and applying the Lefschetz hyperplane theorem for complete
intersections (see for instance [D]), we have $\pi_1(N)\imic\pi(N\cap
H_0\cap H_1)$. But $N\cap H_0\cap H_1$ is given by
$u_0=u_1=u_2+u_3+u_4+u_5=0$ and so is simply-connected. $N$~is
singular, but a resolution $\psi:Y\to N$ is described precisely in
[BN], Section~9, especially (9.1) and (9.3). The morphism $\psi$
contracts only rational varieties (Cayley nodal cubic surfaces,
quadrics and rational curves) and in particular has simply-connected
fibres, so it does not affect the fundamental group.~\qed

Recently Beauville has constructed an example of a Calabi-Yau
threefold whose fundamental group is not abelian ([Be]). Perhaps some
Siegel modular threefolds also have those properties. Unfortunately it
is not easy to tell when a Siegel modular threefold is a Calabi-Yau,
and all the examples given in this paper with non-abelian fundamental
group are of general type.

\beginsection References

\noindent[SC] A. Ash, D Mumford, M. Rapoport \& Y. Tai, {\it Smooth
Compactifications of Locally Symmetric Varieties,} Math. Sci. Press,
Brookline, Mass., 1975.

\noindent[BN] W. Barth \& I. Nieto, Abelian surfaces of type $(1,3)$ and
quartic surfaces with $16$~skew lines, J. Alg. Geom. {\bf 3} (1994)
173-222.

\noindent[Be] A.Beauville, A Calabi-Yau threefold with non-abelian fundamental
group. Preprint 1995.

\noindent[BG] J.W. Bruce \& P. Giblin, {\it Curves and Singularities,}
Cambridge University Press, Cambridge 1992.

\noindent[D] A. Dimca, {\it Singularities and Toplogy of Hypersurfaces,}
Springer, New York 1992.

\noindent[F] W. Fulton, {\it Introduction to toric varieties,}
Princeton University Press, Princeton 1993.

\noindent[vdG] G. van der Geer, {\it Hilbert Modular Surfaces,}
Springer, Heidelberg 1987.

\noindent[G] J. Grosche, \"Uber die Fundamentalgruppe von Quotientr\"aumen
Siegelscher Modulgruppen, J. reine u. angew. Math. {\bf 281} (1976), 53--79.

\noindent[HK] H. Heidrich \& F.W. Kn\"oller, \"Uber die Fundamentalgruppen von
Siegelscher Modulvariet\"aten vom Grade~$2$, Manuscr. Math. {\bf 57} (1987),
249--262.

\noindent[HKW] K. Hulek, C. Kahn \& S.H. Weintraub, {\it Moduli spaces of
abelian surfaces: compactification, degenerations and theta
functions,} De Gruyter, Berlin 1993.

\noindent[HS] K. Hulek \& G.K. Sankaran, The fundamental group of some Siegel
modular threefolds. In {\it Abelian Varieties,} (W. Barth, K. Hulek \&
H. Lange, eds.) De Gruyter, Berlin 1995.

\noindent[K] F.W. Kn\"oller, De Fundamentalgruppen der Siegelscher
Modulvariet\"aten, Abh. Math. Sem. Univ. Hamburg {\bf 57} (1987),
203--213.

\noindent[Na] I. Naruki, On smooth quartic embedding of Kummer surfaces, Proc
Japan Acad. Ser.~A {\bf 67} (1991), 223--225.

\noindent[Y] T. Yamazaki, On Siegel modular forms of degree two, A.J.M. {\bf
98} (1976), 39--53.
\raggedright
\medskip
\noindent G.~K. Sankaran, Department of Pure Mathematics and
Mathematical Statistics,\hfil\break
16,~Mill~Lane, Cambridge CB2~1SB, England.

\end